\input amstex

\catcode`@=11
\nopagenumbers
\vbadness=10000
\hbadness=10000
\hsize=12.5cm
\vsize=19.2cm
\topskip=12pt
\parindent=0.5cm
\newskip\litemindent
\litemindent= 0,7 cm  
\parskip=0pt
\widowpenalty=10000
\clubpenalty=10000
\hfuzz=1.5pt
\abovedisplayskip=6pt plus 1pt
\abovedisplayshortskip=0pt plus 1pt
\belowdisplayskip=6pt plus 1pt 
\belowdisplayshortskip=6pt plus 1pt
\frenchspacing



\font\authorfont=cmti10 at 14.40pt

\font\seventeenrm=cmbx10 at 17.28pt
\font\seventeenit=cmbxti10 at 17.28pt
\font\seventeeni=cmmib10 at 17.28pt
\font\seventeensy=cmbsy10 at 17.28pt
\font\seventeenex=cmex10 at 17.28pt
\font\seventeenmsa=msam10 at 17.28pt
\font\seventeenmsb=msbm10 at 17.28pt
\font\seventeeneuf=eufm10 at 17.28pt

\font\fourteenrm=cmbx10 at 14.40pt
\font\fourteenit=cmbxti10 at 14.40pt
\font\fourteeni=cmmib10 at 14.40pt
\font\fourteensy=cmbsy10 at 14.40pt
\font\fourteenex=cmex10 at 14.40pt
\font\fourteenmsa=msam10 at 14.40pt
\font\fourteenmsb=msbm10 at 14.40pt
\font\fourteeneuf=eufm10 at 14.40pt

\font\twelverm=cmbx10 at 12pt
\font\twelveit=cmbxti10 at 12pt
\font\twelvei=cmmib10 at 12pt
\font\twelvesy=cmbsy10 at 12pt
\font\twelveex=cmex10 at 12pt
\font\twelvemsa=msam10 at 12pt
\font\twelvemsb=msbm10 at 12pt
\font\twelveeuf=eufm10 at 12pt

\font\tenrm=cmr10
\font\tenit=cmti10
\font\tenbf=cmbx10
\font\tenib=cmmib10 
\font\teni=cmmi10
\font\tensy=cmsy10
\font\tenbsy=cmbsy10
\font\tenex=cmex10 
\font\tenmsa=msam10 
\font\tenmsb=msbm10 
\font\teneuf=eufm10

\font\ninerm=cmr9
\font\nineit=cmti9
\font\ninebf=cmbx9
\font\ninei=cmmi9
\font\ninesy=cmsy9

\font\ninemsa=msam9 
\font\ninemsb=msbm10 at 9pt
\font\nineeuf=eufm9

\font\eightrm=cmr8
\font\eightit=cmti8
\font\eightbf=cmbx8
\font\eighti=cmmi8
\font\eightib=cmmib8
\font\eightsy=cmsy8
\font\eightbsy=cmbsy8

\font\eightmsa=msam8
\font\eightmsb=msbm8
\font\eighteuf=eufm8

\font\sevenrm=cmr7

\font\sevenbf=cmbx7
\font\seveni=cmmi7

\font\sevensy=cmsy7

\font\sevenmsa=msam7
\font\sevenmsb=msbm7
\font\seveneuf=eufm7

\font\fiverm=cmr5
\font\fivebf=cmbx5
\font\fivesy=cmsy5

\font\fivei=cmmi5

\font\fivemsa=msam5
\font\fivemsb=msbm5
\font\fiveeuf=eufm5

\newfam\msafam
\newfam\msbfam
\newfam\euffam

\def\seventeenpoint{\def\rm{\fam0\seventeenrm}%
\textfont0=\seventeenrm\scriptfont0=\fourteenrm 
\scriptscriptfont0=\twelverm
\textfont1=\seventeeni\scriptfont1=\fourteeni 
\scriptscriptfont1=\twelvei
\def\mit{\fam1\seventeeni}\def\oldstyle{\fam1\seventeeni}%
\textfont2=\seventeensy\scriptfont2=\fourteensy 
\scriptscriptfont2=\twelvesy
\def\cal{\fam2\seventeensy}%
\textfont3=\seventeenex\scriptfont3=\fourteenex 
\scriptscriptfont3=\twelveex
\def\it{\fam\itfam\seventeenit}%
\textfont\itfam=\seventeenit
\def\bf{\rm}
\def\msa{\fam\msafam\seventeenmsa}%
\textfont\msafam=\seventeenmsa\scriptfont\msafam=\fourteenmsa 
\scriptscriptfont\msafam=\twelvemsa
\def\msb{\fam\msbfam\seventeenmsb}%
\textfont\msbfam=\seventeenmsb\scriptfont\msbfam=\fourteenmsb 
\scriptscriptfont\msbfam=\twelvemsb
\def\euf{\fam\euffam\seventeeneuf}%
\textfont\euffam=\seventeeneuf\scriptfont\euffam=\fourteeneuf 
\scriptscriptfont\euffam=\twelveeuf
\setbox\strutbox=\hbox{\vrule height16pt depth4pt width0pt}%
\baselineskip=20pt\seventeenrm}

\def\fourteenpoint{\def\rm{\fam0\fourteenrm}%
\textfont0=\fourteenrm\scriptfont0=\twelverm 
\scriptscriptfont0=\tenbf
\textfont1=\fourteeni\scriptfont1=\twelvei 
\scriptscriptfont1=\tenib
\def\mit{\fam1\fourteeni}\def\oldstyle{\fam1\fourteeni}%
\textfont2=\fourteensy\scriptfont2=\twelvesy 
\scriptscriptfont2=\tenbsy
\def\cal{\fam2\fourteensy}%
\textfont3=\fourteenex\scriptfont3=\twelveex 
\scriptscriptfont3=\tenex
\def\it{\fam\itfam\fourteenit}%
\textfont\itfam=\fourteenit
\def\bf{\rm}
\def\msa{\fam\msafam\fourteenmsa}%
\textfont\msafam=\fourteenmsa\scriptfont\msafam=\twelvemsa 
\scriptscriptfont\msafam=\tenmsa
\def\msb{\fam\msbfam\fourteenmsb}%
\textfont\msbfam=\fourteenmsb\scriptfont\msbfam=\twelvemsb 
\scriptscriptfont\msbfam=\tenmsb
\def\euf{\fam\euffam\fourteeneuf}%
\textfont\euffam=\fourteeneuf\scriptfont\euffam=\twelveeuf 
\scriptscriptfont\euffam=\teneuf
\setbox\strutbox=\hbox{\vrule height13pt depth4pt width0pt}%
\baselineskip=16pt\fourteenrm}

\def\twelvepoint{\def\rm{\fam0\twelverm}%
\textfont0=\twelverm\scriptfont0=\tenbf 
\scriptscriptfont0=\eightbf
\textfont1=\twelvei\scriptfont1=\tenib 
\scriptscriptfont1=\eightib
\def\mit{\fam1\twelvei}\def\oldstyle{\fam1\twelvei}%
\textfont2=\twelvesy\scriptfont2=\tenbsy 
\scriptscriptfont2=\eightbsy
\def\cal{\fam2\twelvesy}%
\textfont3=\twelveex\scriptfont3=\twelveex 
\scriptscriptfont3=\twelveex
\def\it{\fam\itfam\twelveit}%
\textfont\itfam=\twelveit
\def\bf{\rm}%
\def\msa{\fam\msafam\twelvemsa}%
\textfont\msafam=\twelvemsa \scriptfont\msafam=\tenmsa
\scriptscriptfont\msafam=\eightmsa%
\def\msb{\fam\msbfam\twelvemsb}%
\textfont\msbfam=\twelvemsb \scriptfont\msbfam=\tenmsb%
\scriptscriptfont\msbfam=\eightmsb%
\def\euf{\fam\euffam\twelveeuf}%
\textfont\euffam=\twelveeuf \scriptfont\euffam=\teneuf%
\scriptscriptfont\euffam=\eighteuf%
\setbox\strutbox=\hbox{\vrule height10pt depth4pt width0pt}%
\baselineskip=14pt\rm}

\def\tenpoint{\def\rm{\fam0\tenrm}%
\textfont0=\tenrm\scriptfont0=\sevenrm 
\scriptscriptfont0=\fiverm
\textfont1=\teni\scriptfont1=\seveni 
\scriptscriptfont1=\fivei
\def\mit{\fam1\teni}\def\oldstyle{\fam1\teni}%
\textfont2=\tensy\scriptfont2=\sevensy 
\scriptscriptfont2=\fivesy
\def\cal{\fam2\tensy}%
\textfont3=\tenex\scriptfont3=\tenex 
\scriptscriptfont3=\tenex
\def\it{\fam\itfam\tenit}%
\textfont\itfam=\tenit
\def\bf{\fam\bffam\tenbf}%
\textfont\bffam=\tenbf\scriptfont\bffam=\sevenbf
\scriptscriptfont\bffam=\fivebf
\def\msa{\fam\msafam\tenmsa}%
\textfont\msafam=\tenmsa \scriptfont\msafam=\sevenmsa
\scriptscriptfont\msafam=\fivemsa
\def\msb{\fam\msbfam\tenmsb}%
\textfont\msbfam=\tenmsb \scriptfont\msbfam=\sevenmsb
\scriptscriptfont\msbfam=\fivemsb
\def\euf{\fam\euffam\teneuf}%
\textfont\euffam=\teneuf \scriptfont\euffam=\seveneuf
\scriptscriptfont\euffam=\fiveeuf
\aline=12pt plus 1pt minus 1pt
\halfaline=6pt plus 1pt minus 1pt
\setbox\strutbox=\hbox{\vrule height8.5pt depth3.5pt width0pt}%
\baselineskip=12pt\rm}

\def\ninepoint{\def\rm{\fam0\ninerm}%
\textfont0=\ninerm \scriptfont0=\sevenrm \scriptscriptfont0=\fiverm
\textfont1=\ninei\scriptfont1=\seveni \scriptscriptfont1=\fivei
\def\mit{\fam1\ninei}\def\oldstyle{\fam1\ninei}%
\textfont2=\ninesy \scriptfont2=\sevensy \scriptscriptfont2=\fivesy
\def\cal{\fam2\ninesy}%
\textfont3=\tenex \scriptfont3=\tenex \scriptscriptfont3=\tenex
\def\it{\fam\itfam\nineit}%
\textfont\itfam=\nineit
\def\bf{\fam\bffam\ninebf}%
\textfont\bffam=\ninebf \scriptfont\bffam=\sevenbf
\scriptscriptfont\bffam=\fivebf
\def\msa{\fam\msafam\ninemsa}%
\textfont\msafam=\ninemsa \scriptfont\msafam=\sevenmsa
\scriptscriptfont\msafam=\fivemsa
\def\msb{\fam\msbfam\ninemsb}%
\textfont\msbfam=\ninemsb \scriptfont\msbfam=\sevenmsb
\scriptscriptfont\msbfam=\fivemsb
\def\euf{\fam\euffam\nineeuf}%
\textfont\euffam=\nineeuf \scriptfont\euffam=\seveneuf
\scriptscriptfont\euffam=\fiveeuf%
\aline=11pt plus 1pt minus 1pt
\halfaline=5pt plus 1pt minus 1pt
\setbox\strutbox=\hbox{\vrule height7pt depth3pt width0pt}%
\baselineskip=11pt\rm}

\def\eightpoint{\def\rm{\fam0\eightrm}%
\textfont0=\eightrm \scriptfont0=\sevenrm \scriptscriptfont0=\fiverm
\textfont1=\eighti\scriptfont1=\seveni \scriptscriptfont1=\fivei
\def\mit{\fam1\eighti}\def\oldstyle{\fam1\eighti}%
\textfont2=\eightsy \scriptfont2=\sevensy \scriptscriptfont2=\fivesy
\def\cal{\fam2\eightsy}%
\textfont3=\tenex \scriptfont3=\tenex \scriptscriptfont3=\tenex
\def\it{\fam\itfam\eightit}%
\textfont\itfam=\eightit
\def\bf{\fam\bffam\eightbf}%
\textfont\bffam=\eightbf \scriptfont\bffam=\sevenbf
\scriptscriptfont\bffam=\fivebf
\def\msa{\fam\msafam\eightmsa}%
\textfont\msafam=\eightmsa \scriptfont\msafam=\sevenmsa
\scriptscriptfont\msafam=\fivemsa
\def\msb{\fam\msbfam\eightmsb}%
\textfont\msbfam=\eightmsb \scriptfont\msbfam=\sevenmsb
\scriptscriptfont\msbfam=\fivemsb
\def\euf{\fam\euffam\eighteuf}%
\textfont\euffam=\eighteuf \scriptfont\euffam=\seveneuf
\scriptscriptfont\euffam=\fiveeuf%
\aline=10pt plus 1pt minus 1pt
\halfaline=5pt plus 1pt minus 1pt
\setbox\strutbox=\hbox{\vrule height7pt depth3pt width0pt}%
\baselineskip=10pt\rm}

\skewchar\teni='177
\skewchar\ninei='177 
\skewchar\eighti='177 
\skewchar\seveni='177 
\skewchar\fivei='177
\skewchar\tensy='60
\skewchar\ninesy='60 
\skewchar\eightsy='60 
\skewchar\sevensy='60 
\skewchar\fivesy='60


\long\def\Title#1\par{%
\global\titlepagetrue
{\parindent=0pt
\raggedcenter\pretolerance=10000
\seventeenpoint #1\par}
\vskip24pt}

\long\def\Author#1\par{%
\centerline{\fourteenpoint\authorfont #1}
\vskip60pt
\tenpoint\noindent\ignorereturn}

\def\Section{\removelastskip
\goodbreak\vskip36pt plus 1pt minus 6pt \section}

\def\section#1\par{%
{\raggedcenter
\interlinepenalty=10000\pretolerance=10000
\noindent\fourteenpoint #1\nobreak\par\nobreak}
\nobreak\vskip12pt\nobreak
\noindent\tenpoint\rm\ignorereturn}

\def\subsection#1\par{%
{\raggedcenter
\interlinepenalty=10000\pretolerance=10000
\noindent\twelvepoint #1\nobreak\par\nobreak}
\nobreak\vskip12pt\nobreak
\noindent\tenpoint\rm\ignorereturn}

\def\References #1 {\ifdim\lastskip<25pt \removelastskip
\vskip24pt plus 1pt\fi
\setbox0=\hbox{\ninepoint [#1]\enspace}%
\centerline{\twelverm References}\nobreak\par\nobreak
\interlinepenalty=10000\parskip=0pt plus 1pt
\litemindent=\wd0
\ninepoint
\nobreak\vskip7pt\nobreak}



\newbox\footbox
\setbox\footbox=\hbox{\ninerm 22)~}

\newdimen\footnotespace
\footnotespace=0pt

\def\footnote#1{\let\@sf\empty 
  \ifhmode\edef\@sf{\spacefactor\the\spacefactor}\/\fi
$^{#1}$\@sf\vfootnote{#1}}

\def\vfootnote#1{\insert\footins\bgroup
\interlinepenalty\interfootnotelinepenalty
\splittopskip\ht\strutbox 
\splitmaxdepth\dp\strutbox \floatingpenalty\@MM
\leftskip\z@skip \rightskip\z@skip \spaceskip\z@skip \xspaceskip\z@skip\parindent=\wd\footbox
\litem{\ninerm #1}\ninepoint\footstrut\futurelet\next\fo@t}
\def\fo@t{\ifcat\bgroup\noexpand\next \let\next\f@@t
\else\let\next\f@t\fi \next}
\def\f@@t{\bgroup\aftergroup\@foot\let\next}
\def\f@t#1{#1\@foot}
\def\@foot{\strut\egroup}
\def\footstrut{\vbox to\splittopskip{}}
\skip\footins=\bigskipamount 
\count\footins=1000 
\dimen\footins=8in 


\newif\iftitlepage

\def\makeheadline{\iftitlepage \global\titlepagefalse
\vbox{\titleheadline}%
\vskip12pt
\else \vbox{\ifodd\pageno\rightheadline\else\leftheadline\fi}%
\vskip12pt\fi}

\def\lefthead{}
\def\righthead{}

\def\rightheadline{\line{\vbox to 8.5pt{}%
\hphantom{\tenrm\folio}\ninepoint\hfill\righthead\hfill{\tenrm\folio}}}
\def\leftheadline{\line{\vbox to 8.5pt{}%
\ninepoint{\tenrm\folio}\hfill\lefthead\hfill\hphantom{\tenrm\folio}}}
\def\titleheadline{\line{\vbox to 8.5pt{}%
\ninepoint\it\hfill[Page \folio]}}

\output{\plainoutput}


\def\Litem#1#2{\par\noindent\hangindent#1\litemindent
\hbox to #1\litemindent{\hfill\hbox to \litemindent
{#2 \hfill}}\ignorespaces}
\def\litem{\Litem1}

\newskip\aline \newskip\halfaline
\aline=12pt plus 1pt minus 1pt
\halfaline=6pt plus 1pt minus 1pt
\def\skipaline{\vskip\aline}

\def\ignorereturn{\def\neext{\afterassignment\restart
       \let\next}\neext}
\def\restart{\ifx\next\par\else\let\neext\next\fi\neext}

\def\raggedcenter{\leftskip=0pt plus 4em \rightskip
=\leftskip \parfillskip=0pt \spaceskip=.3333em
\xspaceskip=.5em \pretolerance=9999 \tolerance=9999
\hyphenpenalty=9999 \exhyphenpenalty=9999 }

\def\dotfill{\leaders\hbox to 1em{\hss.\hss}\hfill}

\def\Classification{\ifdim\lastskip<\aline\removelastskip\skipaline\fi
\noindent 1991 Mathematics Subject Classification: }

\def\Abstract{\ninepoint\noindent\bf Abstract. \rm}

\long\def\Proclamation#1. #2\par{\ifdim\lastskip<\aline\removelastskip
\penalty-250 \skipaline\fi{\def\\##1){\litem{\rm(##1)}}\noindent
\bf#1\unskip. \it#2\par}\skipaline}

\def\Theorem {\Proclamation Theorem }

\def\Corollary {\Proclamation Corollary }
\def\Lemma {\Proclamation Lemma }

\def\Proof{\ifdim\lastskip<\aline\removelastskip\skipaline\fi
\noindent\it Proof. \rm}

\def\qedbox{$\rlap{$\sqcap$}\sqcup$}
\def\qed{\nobreak\hfill\penalty250 \hbox{}\nobreak\hfill\qedbox\skipaline}

\def\proclamation#1. {\ifdim\lastskip<\aline\removelastskip\penalty-250
      \skipaline\fi\noindent\bf#1\unskip. \rm}

\def\ref #1 {\vskip4pt\litem{[#1]}}

\catcode`@=12
\def\@{\hbox{-}}

\catcode`@=\active

\UseAMSsymbols
%
\redefine\le{\leqslant}
\redefine\ge{\geqslant}
\define\a{\alpha}

\define\({\left(}
\define\){\right)}
\define\pfrac#1#2{\( \frac{#1}{#2} \)}

\Title On two conjectures of Sierpi\'nski concerning the arithmetic functions
$\sigma$ and $\phi$

{\fourteenpoint\authorfont
\centerline{Kevin Ford}
\centerline{Sergei Konyagin$^*$}}
\vskip 30pt
\vfootnote{*}{The second author was supported by NSF grant DMS 9304580.}

{\tenpoint \it \noindent
 Dedicated to Professor Andrzej Schinzel on the occasion of
 his 60th birthday. }

\vskip 30pt
\Abstract  Let $\sigma(n)$ denote the sum of the positive divisors of $n$.
In this note it is shown that for any positive integer $k$, there is
a number $m$ for which the equation $\sigma(x)=m$ has exactly $k$ solutions,
settling a conjecture of Sierpi\'nski.   Additionally, it is shown that
for every positive even $k$, there is a number $m$ for which the equation
$\phi(x)=m$ has exactly $k$ solutions, where $\phi$ is Euler's function.

\Classification Primary 11A25, 11N64; Secondary 11N36.

\Section 1. Introduction

\noindent
For each natural number $m$, let $A(m)$ denote the number of solutions  of
$\phi(x)=m$ and let $B(m)$ denote the number of solutions of $\sigma(x)=m$.
Here $\phi(x)$ is Euler's function and $\sigma(x)$ is the sum of divisors
function.  About 40 years ago, Sierpi\'nski made two conjectures about the
possible values of $A(m)$ and $B(m)$ (see [S1],
[E,p. 12] and Conjectures $C_{14}$ and $C_{15}$ of [S2]).

\Proclamation Conjecture 1 (Sierpi\'nski).
For each $k\ge 2$, there is a number $m$ with $A(m)=k$.

\Proclamation Conjecture 2 (Sierpi\'nski).
For each $k\ge 1$, there is a number $m$ with $B(m)=k$.

An older conjecture of
Carmichael [C1,C2] states that $A(m)$ can never equal 1.
Carmichael's Conjecture remains unproven, however it is known that a
counterexample $m$ must exceed $10^{10^{10}}$ (c.f. Theorem 6 and
section 7 of [F1]).

Both of Sierpi\'nski's conjectures were deduced by Schinzel
[S1] as a consequence of his Hypothesis H [SS].

\Proclamation {Schinzel's Hypothesis H}.
Suppose $f_1(n), \ldots, f_k(n)$ are irreducible, integer valued polynomials
(for integral $n$) with positive leading coefficients.  Also suppose that
for every integer $q\ge 2$, there is an integer
 $n$ for which $q$ does not divide $f_1(n)\cdots f_k(n)$.
Then the numbers $f_1(n),\ldots,f_k(n)$ are simultaneously prime
for infinitely many positive integers $n$.

By an inductive approach, the first author [F1,Lemma 7.1] has shown that
Conjectures 1 and 2 follow from Dickson's Prime $k$-tuples Conjecture [D],
which is the special case of Hypothesis H when each $f_i(n)$ is linear.

Although Hypothesis H has not been proved in even the simplest case of
two linear polynomials (generalized twin primes), sieve methods have shown the
conclusion to hold if the numbers $f_1(n),\ldots,f_k(n)$ are allowed to be
primes or ``almost primes'' (non-primes with few prime factors).  See [HR] 
for specifics.
Taking a new approach we utilize these almost primes to
prove Conjecture 2 unconditionally.
The same method is applicable to Conjecture 1, but
falls short of a complete proof because of the
(probable) non-existence
of a number with $A(m)=1$.  The fact that $B(1)=1$ is
crucial to the proof of Conjecture 2. 

\Theorem 1.  For every $k\ge 1$, there is a number $m$ with $B(m)=k$.

\Theorem 2.  Suppose $r$ is a positive integer and $A(m)=k$.  Then
there is a number $l$ for which $A(lm)=rk$.

\Corollary 3.  If $A(m)=k$ is known to be solvable for $2\le k\le C$,
then $A(m)=k$ has a solution for every
$k$ divisible by a prime $\le C$.  In particular, $A(m)=k$ is solvable
for all even $k$.

The first author has succeeded in proving
Conjecture 1 for all $k\ge 2$ by combining the inductive 
approach in [F1] with the theory of almost primes.
The details are very complex and will appear in a forthcoming paper
[F2].

\Section 2. Preliminary lemmas

\noindent
Let $\omega(n)$ denote the number of distinct prime factors of $n$, let
$P^-(n)$ denote the smallest prime factor of $n$, and let $[x]$ denote the
greatest integer $\le x$.
The first two lemmas provide the construction of numbers $m$ with a desired
value of $A(m)$ or $B(m)$.

\Lemma 1.
Suppose $A(m)=k$, $r\ge 2$, $n\ge 2$ and $p_{i,j}$ $(i=1,\ldots,r;j=1,\ldots,
n)$ are primes larger than $2^rm+1$.  For each $i$, let $q_i=p_{i,2}p_{i,3}
\cdots p_{i,n}$, and let $t$ be the product of all primes $p_{i,j}$.
Suppose further that
\\i) $2p_{i,1}q_j+1$ is prime whenever $i=1$, $j=1$ or $j=i$,
\\ii) no $p_{i,j}$ equals any of the primes listed in (i),
\\iii) except for the numbers listed in (i), for each $d_1|t$
with $d_1>1$ and $d_2|2^{r-1}m$, $2d_1d_2+1$ is composite.
\endgraf\noindent Then $A(2^rtm)=rk$.

\Proof Suppose that
$\phi(x)=2^rtm$.  No $p_{i,j}$ may divide $x$, for otherwise
$p_{i,j}-1|2^rtm$, which is impossible by conditions (ii),
(iii) and the fact that each $p_{i,j}>2^rm+1$.
Therefore, each $p_{i,j}$ divides a number $s_{i,j}-1$,
where $s_{i,j}$ is a prime divisor of $x$. 
Therefore, $s_{i,j}= d p_{i,j}+1$, where $d|2^rmt/p_{i,j}$ and $2|d$.  By
condition (iii), $s_{i,j}$ must be one of the primes listed in (i) and
by condition (ii), each prime $s_{i,j}$ divides $x$ to the first power only.
By (i), there are $r$ choices for $s_{1,1}$ and once $s_{1,1}$ is chosen
the other primes $s_{i,j}$ are uniquely determined.   For each choice,
$$
\phi(s_{1,1}s_{2,1}\cdots s_{r,1}) = 2^rt,
$$
and thus $\phi(x/(s_{1,1}\cdots s_{r,1}))=m$, which has exactly $k$ solutions.
\qed

\Lemma 2.
Suppose $r\ge 2$, $n\ge 2$ and $p_{i,j}$ $(i=1,\ldots,r;j=1,\ldots,
n)$ are primes larger than $2^r+1$.  For each $i$, let $q_i=p_{i,2}p_{i,3}
\cdots p_{i,n}$, and let $t$ be the product of all primes $p_{i,j}$.
Suppose further that
\\i) $2p_iq_j-1$ is prime whenever $i=1$, $j=1$ or $j=i$,
\\ii) $\sigma(\pi^b)\nmid 2^rt$ for every prime $\pi$
and integer $b\ge 2$ with $\sigma(\pi^b)>2^r$,
\\iii) except for the numbers listed in (i), for each $d_1|t$
with $d_1>1$ and $d_2|2^{r-1}$, $2d_1d_2-1$ is composite.
\endgraf\noindent Then $B(2^rt)=r$.

\Proof Suppose that $\sigma(x)=2^rt$.
Each $p_{i,j}$ divides a number $\sigma(s_{i,j}^b)$,
where $s_{i,j}^b$ is a prime power divisor of $x$.
Condition (ii) implies $b=1$, so 
$s_{i,j}= d p_{i,j}-1$, where $d$ is an even divisor of
 $2^rt/p_{i,j}$.  By
condition (iii), $s_{i,j}$ must be one of the primes listed in (i).
There are $r$ choices for $s_{1,1}$ and once $s_{1,1}$ is chosen
the other primes $s_{i,j}$ are uniquely determined.   For each choice,
$$
\sigma(s_{1,1}s_{2,1}\cdots s_{r,1}) = 2^rt,
$$
which forces $x=s_{1,1}\cdots s_{r,1}$.
\qed

To show such sets of primes $(p_{i,j})$ exist, the first tool we
require is a lower bound on the density of primes $s$ for which
$\frac{s-1}2$ (or $\frac{s+1}2$) is an almost prime.

\Lemma 3.
Let $a=1$ or $a=-1$.  For some positive $\a$ and $x$ sufficiently large,
there are $\gg x/\log^2 x$ primes $x/2< s\le x$ for which $s=2u+a$,
$u$ has at least 2 prime factors and every prime factor of $u$ exceeds
$x^\a$.

\Proof
This follows from the linear sieve 
and the Bombieri-Vinogradov prime number theorem (Lemma 3.3 of [HR]) to
bound the error terms.  By Theorem 8.4 of [HR], we have
$$
\# \{x/2< s\le x: s, \tfrac12(s-a) \text{ both prime} \} \le
(4+o(1))\frac{x}{\log^2 x}
$$
and for $x\ge x_0(\a)$
$$
\# \{x/2<s\le x: s\text{ prime }, P^-(\tfrac12(s-a)) >x^\a \} \ge
\( \frac{e^{-\gamma}}{\a} f(1/(2\a)) + o(1) \)\frac{x}{\log^2 x},
$$
where $f$ is the usual lower bound sieve function and $\gamma$ is the
Euler-Mascheroni constant.  Taking $\a=\frac18$ and noting that
$f(4)=\frac12 e^{\gamma} \log 3$, the number of primes $x/2<s\le x$
for which $u=\frac12(s-a)$ contains at least 2 prime factors and all
prime factors of $u$ exceed $x^\a$ is at least $0.39 x/\log^2 x$ for
large $x$.
\qed

In the argument below it is critical
that the numbers $\frac12 (s-a)$ have at least two prime factors.
This may be the first application of lower
bound sieve results where almost primes are desired and primes are not.

\Lemma 4.
Suppose $g\ge 1$, and $a_i,b_i (i=1,\ldots,g)$ are integers
satisfying
$$
E := \prod_{i=1}^g a_i \prod_{1\le r<s\le g} (a_rb_s-a_sb_r) \ne 0.
$$
Let $\rho(p)$ denote the number of solutions of
$$
\prod_{i=1}^g (a_in+b_i) \equiv 0\pmod{p},
$$
and suppose $\rho(p)<p$ for every prime $p$.  If $\log E \ll \log z$, then
the number of $n$ with $z<n\le 2z$ and $P^-(a_in+b_i)>z^\a$
for $i=1,\ldots,g$ is
$$
\split
&\ll_{g,\a} \frac{z}{\log^g z} \prod_p \( 1 - \frac{\rho(p)-1}{p-1} \)
\(1-\frac1{p} \)^{1-g} \\
&\ll_{g,\a} \frac{z}{\log^g z} \( \frac{E}{\phi(E)} \)^g
\ll_{g,\a} \frac{z(\log\log z)^g}{\log^g z}.
\endsplit
$$

\Proof  This is essentially Theorem 5.7 of [HR].  The second
part follows from the fact that $\rho(p)=g$ unless $p|E$, in which case
$\rho(p)<g$.
\qed

\Lemma 5. For any real $\beta>0$,
$$
\sum_{k\le x} \pfrac{k}{\phi(k)}^\beta \ll_\beta x.
$$

\Proof Write $(k/\phi(k))^\beta=\sum_{d|k} g(d)$, where $g$ is the
multiplicative function satisfying $g(p)=(p/(p-1))^\beta-1$ for primes $p$
and $g(p^a)=0$ when $a\ge 2$.  Then
$$
\sum_{k\le x} (k/\phi(k))^\beta = \sum_{d\le x} g(d)[x/d] \le x\prod_p
(1+g(p)/p) = c(\beta)x.
$$

\Section 3.  The main argument

Fix $a=1$ or $a=-1$.  The primes $s$ counted in
Lemma 3 have the property that $\omega(\frac12(s-a))\le
[1/\a]$.
Therefore, there exists a number $n$ $(1\le n\le [1/\a]-1)$ and
some pair $y,z$ with $x/16 \le yz \le x/2$, $y>x^\a$
such that
$$
\#\{y<p\le 2y, z<q\le 2z: p,2pq+a\text{ prime}, \omega(q)=n,
P^-(q)>y \} \gg \frac{x}{\log^3 x}.
$$
Denote by $B$ the set of such pairs ($p,q$).
From now on variables $p,p_i$ will denote primes in $(y,2y]$ and variables
$q,q_i$ will denote numbers in $(z,2z]$ with $n$ prime factors, each exceeding
$y$.  Implied constants
in the following may depend on $r$, $n$ or $m$.

\Lemma 6.  The number of $2r$-tuples $(p_1,\ldots,q_r)$
with each $(p_i,q_i)\in B$
which satisfy condition (i) but fail condition (ii) or (iii) (referring
either to Lemma 1 or Lemma 2 and writing $p_i=p_{1,i}$ and
$q_i=p_{i,2}\cdots p_{i,n}$) is
$$
\ll \frac{x^{r}(\log\log x)^{rn+4r-1}}{(\log x)^{5r-1}}.
$$

\Proof
We first count those $2r$-tuples satisfying (i) but failing (ii).
When $a=1$, all of the $2r$-tuples satisfy condition (ii) in Lemma 1, since 
$2p_{i,1}q_j+1 \gg x$ and each $p_{i,j} \ll x^{1-\a}$.
If condition (ii) of Lemma 2 fails, then 
$y/2 \le \pi^b \le 2^r t \le (2x)^r$.  Therefore,
the number of $2r$-tuples
not satisfying (ii) is bounded above by
$$
\sum_{y/2 \le \pi^b \le (2x)^r} \frac{(2x)^r}{\pi^b} 
\ll x^r \sum_{b=2}^{\infty} \sum_{\pi \ge (y/2)^{1/b}} \frac{1}{\pi^b}
\ll x^{r-\a/2}.
$$
Counting the $2r$-tuples satisfying (i) but failing (iii) is a
straightforward application of Lemma 4.  First fix
$d_2$ and the set of pairs ($i,j$) for which
$p_{i,j}|d_1$  (there are finitely many such choices).
Each of the numbers listed in (i) and (iii) are linear in
all the variables $p_{i,j}$, thus applying Lemma 4 successively with the
variables $p_{i,j}$ (in some order) gives the desired upper bound on their
number.

We illustrate this process in the case $r=3$, $n=2$,
$d_1=p_{2,2}p_{2,3}p_{3,3}$, $d_2$ arbitrary.
Fix distinct primes $p_{1,2},p_{1,3},p_{2,2},p_{2,3},p_{3,2}$.  Since
$p_{3,2} \ll z/y$,  by Lemma 4 the number of primes $p_{3,3}$ such that
$2d_2 p_{2,2}p_{2,3}p_{3,3}+a$ is prime is
$$
\ll \frac{z(\log\log x)^2}{p_{3,2} \log^2 x}.
$$
Given $p_{3,3}$ (i.e. $q_1,q_2,q_3$ are fixed), the number of $p_1$ with
$2p_1q_j+a$ prime ($j=1,2,3$) is $O(y(\log\log x)^4/\log^4 x)$, the number
of $p_2$ with $2p_2q_j+a$ prime ($j=1,2$) is 
$O(y(\log\log x)^3/\log^3 x)$ and the number of $p_3$ with $2p_3q_j+a$
prime ($j=1,3$) is $O(y(\log\log x)^3/\log^3 x)$.  Multiplying these
together and summing over all $p_{i,j}$ ($i=1,2,3; j=2,3$) gives an upper
bound of $O(x^3 (\log\log x)^{13}/\log^{14} x)$ 6-tuples.
\qed

\Lemma 7. The number of $2r$-tuples $(p_1,\ldots,q_r)$, with
each $(p_i,q_i)\in B$,
satisfying condition (i) of Lemma 1 or Lemma 2 is
$$
\gg \frac{x^r}{(\log x)^{5r-2}}.
$$

\Proof
Denote by $P_j$ a
generic $j$-tuple $(p_1, \ldots, p_j)$ with $p_1,\ldots, p_j$ distinct.
Let $N_j(q)$ be the number of $P_j$ such that $2p_iq+a$
is prime for each $i$, and let $M_j(P_j)$ be the number of $q$ such
that $2p_iq+a$ is prime for each $i$.

By the definition of $B$, we have
$$
\sum_q N_1(q) =|B| \gg x/\log^3 x.
$$
Therefore, by H\"older's inequality,
$$
\aligned
S:= \sum_{P_r} M_r(P_r) &= \sum_q N_r(q) = r! \sum_q \binom{N_1(q)}{r} \\
&\gg \sum_{N_1(q)\ge r+1} N_1(q)^r \gg (z/\log z)^{1-r} \(
\sum_{N_1(q)\ge r+1} N_1(q) \)^r \\
&\gg \frac{x^r}{z^{r-1} (\log x)^{2r+1}}.
\endaligned\eqno(1)
$$
Lemma 4 gives
$$
M_r(P_r) \ll L(P_r) \frac{z}{(\log x)^{r+1}}, \eqno(2)
$$
where 
$$
L(P_j) := \prod_{1\le g<h\le j} \frac{|p_g-p_h|}{\phi(|p_g-p_h|)}.\eqno(3)
$$
This follows from the fact that $r+1\ge \rho(p) \ge r+1-k_p$, where $k_p$
is the number of pairs ($i,j$) with $i>j$ and $|p_i-p_j|$ divisible by $p$.
Let $A$ be the number of $p$, so that $A \asymp y/\log x$.
Let $R(k;x)$ denote the number of primes $p\le x-k$ for which $p+k$ is also
prime.  By Lemma 4, when $k\le x/2$ we have
$$
R(k;x) \ll \frac{x}{\log^2 x}\frac{k}{\phi(k)}.
$$
Lemma 5 now gives
$$
\sum_{y<p_1<p_2\le 2y} L(p_1,p_2)^\beta \le
\sum_{k\le y} \pfrac{k}{\phi(k)}^\beta R(k;2y)
\ll_\beta A^2.
$$
Let $H=\binom{j}2$.  Together with (3) and
H\"older's inequality, we have
$$
j! \binom{A}{j} \le \sum_{P_j} L(P_j) \le
 \prod_{1\le g<h\le j} \( A^{j-2} \sum_{p_g,p_h} L(p_g,p_h)^{H} \)^{1/H}
\ll_j A^j\eqno(4)
$$
and similarly
$$
\sum_{P_j} L^2(P_j) \ll_j A_j.
$$
The upper bounds
$$
S \ll \frac{zA^r}{(\log x)^{r+1}}
$$
and
$$
\sum_{P_r} M^2_r(P_r) \ll S^2 A^{-r} \eqno(5)
$$
now follow from  (1), (2) and (4).
Choose $\delta_0>0$ small enough so that 
$$
r! \binom{A}{r}\frac{\delta_0 z} {(\log x)^{r+1}} \le \frac{S}2
$$
and let $P$ denote the set of $P_r$ with
$$
M_r(P_r) \ge \frac{\delta_0 z}{(\log x)^{r+1}}.
$$
By (5) and the Cauchy-Schwarz inequality,
$$
\align
S &\le \(r! \binom{A}{r}-|P|\) \frac{\delta_0}
{(\log x)^{r+1}} + \sum_{P_r\in P} M_r(P_r) \\
&\le \frac{S}2 + O\( |P|^{1/2} S A^{-r/2} \),
\endalign
$$
whence
$$
|P| \gg A^r. \eqno(6)
$$

For each $P_j$, let $J_j(P_j)$ denote
the number of $P_{r-j}$ with $P_{r-j} \cap P_j = \emptyset$ and
$(P_j,P_{r-j}) \in P$.  
Let $\delta_1$ and $\delta_2$ be sufficiently small positive constants,
depending on $r$, but not on $A$.  Let $R$ denote the set of $p$ such that
$J_1(p) \ge \delta_1 A^{r-1}$.  By (6),
if $\delta_1$ is small enough then $|R| \gg A$.
If $p\in R$, denote by $T(p)$ the set of $p'$ such that $J_2(p,p')\ge
\delta_2 A^{r-2}$. 
 If $\delta_2$ is small enough, $|T(p)| \gg A$ uniformly in $p$.
Choose $\delta_2$ so that $\delta_2 < \frac1{2r} \delta_1$.
We first show that
$$
\sum_{\Sb p_1\in R \\ p_2,\ldots,p_r \in T(p_1) \\ (p_1,\ldots,p_r)\in P
\endSb} M_r(p_1,p_2,\ldots, p_r) \gg A^r \frac{z}{(\log x)^{r+1}}. \eqno(7)
$$
The functions $M_j$ are symmetric in all variables, hence
$$
\align
\# \{(p_1,\cdots,p_r)\in P&: p_1\in R; p_2,\ldots,p_r\in T(p_1) \} \\
&\ge \sum_{p_1\in R} J_1(p_1) - r\sum_{\Sb p_1\in R \\ p_2\not\in T(p_1)\endSb}
J_2(p_1,p_2) \\
&\ge |R|\delta_1 A^{r-1} - r|R|A (\delta_2 A^{r-2}) \\
&\ge \tfrac12 |R|\delta_1 A^{r-1} \gg A^r.
\endalign
$$
Together with the definition of $P$, this proves (7).
Next, if $p_1\in R$ and $p_2\in T(p_1)$, then by Lemma 4,
$$
\align
J_2(p_1,p_2) \delta_0 \frac{z}{(\log x)^{r+1}} &\le
\sum_{p_3,\ldots, p_r} M_r(p_1,\ldots,p_r) \\
&= \sum_{q\text{ counted in }M_2(p_1,p_2)} N_{r-2}(q) \\
&\ll \( \frac{y}{\log^2 x} \)^{r-2} M_2(p_1,p_2),
\endalign
$$
whence
$$
M_2(p_1,p_2) \gg \frac{z}{\log^3 x}.\eqno(8)
$$
Let $E$ denote the number of $2r$-tuples $(p_1,\ldots,q_r)$ with
$\gcd(q_i,q_j)>1$ for some $i\ne j$.
Using (7) and (8), the number of $2r$-tuples satisfying condition (i) of
either Lemma 1 or Lemma 2 is at least
$$
\align
 \sum_{p_1,\cdots,p_r}& M_r(p_1,\ldots,p_r) \prod_{j=2}^r M_2(p_1,p_j) -E \\
&\ge \sum_{\Sb p_1\in R \\ p_2,\ldots,p_r\in T(p_1) \endSb} 
M_r(p_1,\ldots,p_r) \prod_{j=2}^r M_2(p_1,p_j) -E \\
&\gg \( \frac{z}{\log^3 x} \) ^{r-1} \sum_{\Sb p_1\in R \\
 p_2,\ldots,p_r\in T(p_1) \endSb}  M_r(p_1,\ldots,p_r) -E\\
&\gg \frac{x^r}{(\log x)^{5r-2}}-E.
\endalign
$$
Trivially $E\ll \frac{x^r}{y}$ and the lemma follows.
\qed

For every $r\ge 2$,
Lemmas 6 and 7 guarantee the existence of
 a set of primes $(p_{i,j})$ satisfying the
hypotheses of Lemma 1 or Lemma 2.  This completes the proof of Theorems 1 and
2.

The methods of this paper also apply to a wide class of multiplicative 
arithmetic functions.  An exposition of some results will appear in section
9 of [F1].


\References HR

\ref C1 Carmichael, R. D., {\it On Euler's $\phi$-function}.
Bull. Amer. Math. Soc. {\bf 13} (1907), \hbox{241--243.}

\ref C2 --- {\it Note on Euler's $\phi$-function}.
Bull. Amer. Math. Soc. {\bf 28} (1922), 109--110.

\ref D Dickson, L. E., {\it A new extension of Dirichlet's
theorem on prime numbers}. Messenger of Math. {\bf 33} (1904),
155--161.

\ref E Erd\"os, P, {\it Some remarks on Euler's $\phi$-function}.
Acta Arith. {\bf 4} (1958), 10--19.

\ref F1 Ford, K, {\it The distribution of totients}.
The Ramanujan J. {\bf 2}, nos. 1--2 (1998), 67--151.

\ref F2 --- {\it The number of solutions of $\phi(x)=m$} Annals of Math.
(to appear)

\ref HR Halberstam, H., Richert, H.-E., Sieve Methods
Academic Press, London 1974.

\ref SS  Schinzel, A.,Sierpi\'nski, W.,
{\it Sur certaines hypoth\`eses concernant les nombres premiers},
Acta Arith. {\bf 4} (1958), 185--208.

\ref S1  Schinzel, A., {\it Sur l'equation $\phi(x)=m$}, 
Elem. Math. {\bf 11} (1956), 75--78.

\ref S2  --- {\it Remarks on the paper ``Sur
certaines hypoth\`eses concernant les nombres premiers''},  Acta Arith.
{\bf 7} (1961/62), 1--8.

\end